\documentclass[12pt, reqno]{amsart} 
\usepackage[top=1in, bottom=1in, left=1in, right=1in]{geometry}

\usepackage{times}
\usepackage[table,dvipsnames]{xcolor}
\PassOptionsToPackage{hyphens}{url}\usepackage[colorlinks=true, linkcolor=black, citecolor=blue, urlcolor=blue]{hyperref}   %% linkable cross-refs
\numberwithin{equation}{section}

\usepackage{scrextend}
\usepackage{amssymb,mathrsfs,amsmath}
\usepackage{mathtools}
\usepackage{bbm}
\usepackage{url}

\theoremstyle{plain}
\newtheorem{thm}{Theorem}
\newtheorem{prop}{Proposition}[section]
\newtheorem{lem}[prop]{Lemma}

\theoremstyle{definition}

\theoremstyle{remark}
\newtheorem*{rem*}{Remark}
\newtheorem*{rems*}{Remarks}

\newcommand{\N}{\mathbb{N}}

\let\lcm\relax% Set equal to \relax so that LaTeX thinks it's not defined
\DeclareMathOperator{\lcm}{lcm}

\newcommand{\bv}\boldsymbol{}

\newcommand{\eps}\varepsilon

\renewcommand{\leq}{\leqslant}

\renewcommand{\geq}{\geqslant}

\begin{document}

\title{On the intervals for the non-existence of covering systems with distinct moduli}

\author{Jack Dalton}
\address{Department of Mathematics\\
	University of Colorado Boulder\\
	2300 Colorado Ave.\\
	Boulder, CO 80309\\
	USA}
\email{{\tt jack.dalton-1@colorado.edu}}

\author{Nic Jones}
\address{Department of Mathematics\\
	College of Charleston\\
	175 Calhoun St.\\
	Charleston, SC 29401\\
	USA}
\email{{\tt jonesna1@g.cofc.edu}}

%%\keywords{}
\date{\today}

\begin{abstract}
In a research seminar in $2006$, M. Filaseta, O. Trifonov, and G. Yu showed for each integer $n\geq3$ there is no distinct covering with all moduli in the interval $[n, 6n]$. In $2022$, this interval was subsequently improved to $[n, 8n]$ by the first author and O. Trifonov. The first author then improved this bound to $[n, 9n]$ in $2023$. Building off their method, we show that for each integer $n\geq 3$, there does not exist a distinct covering system with all moduli in the interval $[n, 10n]$.
\end{abstract}

\maketitle
\section{Introduction}

A \textit{covering system} (also known as a \textit{covering}) $\mathcal{C}$ is a set of congruences $x\equiv a_i\pmod{n_i},\\ \ i\in\{1, 2,\ldots, k\}$, such that every integer satisfies at least one of the congruences.
In this paper we only concern ourselves with finite sets. Further, we say a \textit{minimal} covering system, denoted $\mathcal{C}_m$, is a covering such that no proper subset of it is a covering.\\
\indent In 1950, P. Erd\H{o}s \cite{Erdos} introduced the study of covering systems of the integers. At the time he was interested in distinct covering systems, that is, a covering such that $1<n_1<n_2<\ldots<n_k$. As a result of his work, we are particularly interested in studying covering systems with all distinct moduli.\\
\indent In the same paper, Erd\H{o}s also conjectured if the minimum modulus of a covering is arbitrarily large. This conjecture of Erd\H{o}s became known as the \textit{minimum modulus problem}. In 2007, M. Filaseta, K. Ford, S. Konyagin, C. Pomerance, and G. Yu \cite{Filaseta1} laid the foundation for solving the \textit{minimum modulus problem}. Building off the their work, R. Hough \cite{Hough} in 2015 achieved a substantial breakthrough and solved the minimum modulus problem. In particular, Hough showed that the minimum modulus of any distinct covering system is bounded above by $10^{16}$. This upper bound on the minimum modulus has since been improved by P. Balister, B. Bollob\'as, R. Morris, J. Sahasrabudhe, and M. Tiba \cite{Balister1} to $616000$. We use this bound in this paper.\\
\indent In 1971, C. Krukenberg \cite{Krukenberg} wrote a Ph.D. dissertation on covering systems, but none of his results were ever published in mathematical journals.
Because of this, we outline some of the necessary results of his dissertation that we use in this paper. Firstly, Krukenberg proved the following theorem,
\begin{thm}[Krukenberg, 1971]\label{thm1}
    (i) If the minimum modulus of a distinct covering system is $3$, then its largest modulus is $\geq36$.\\
    \indent (ii) If the minimum modulus of a distinct covering system is $4$, then its largest modulus is $\geq60$.
\end{thm}
Krukenberg states as above that given $m=4$, the largest modulus is at least $M=60$ but he writes ``this result will not be proved here''. In 2022, it was proven by O. Trifonov and the first author that $M=60$ is indeed least possible when $m=4$ \cite{Dalton2}. In addition, Krukenberg proved that if $m=5$, then $M\geq108$. This value of $M=108$ was conjectured by Krukenberg to be the least possible value of $M$ in this case. The first author remarked in his own Ph.D. dissertation that ``We can show that if the least modulus of a distinct covering system is 5, then its largest modulus is at least 84. However, the result is too weak and the proof too long, to be included in this paper''\ \cite{Dalton1}.\\
\indent Lastly, M. Filaseta, O. Trifonov, and G. Yu \cite{Filaseta2} showed in a research seminar in $2006$ that for each integer $m\geq3$ there is no distinct covering system with all moduli in the interval $[m, 6m]$. This interval was subsequently improved by the first author and O. Trifonov in 2022 to $[m, 8m]$ for $m\geq3$. In 2023, the first author proved the same statement for $[m, 9m]$ for $m\geq3$. We prove this result for $[m, 10m]$ for $m\geq3$.
\begin{thm}\label{thm2}
    For each integer $m\geq3$, there is no distinct covering system with all moduli in the interval $[m, 10m]$.
\end{thm}
\indent We organize this paper as follows. In Section 2, we introduce notation and the corollaries and lemmas of Krukenberg and the first author which we use to prove Theorem \ref{thm2}. In Section 3, we prove Theorem \ref{thm2}. Lastly, in Section 4, we state some open problems and discuss possible extensions to the non-existence of distinct coverings with all moduli in the interval $[n, kn]$ for $k>10$.

\section{Reduction of covering systems}
In this section, we introduce the idea of ``reducing'' covering systems. In particular, given a covering $\mathcal{C}$, we can reduce the number of congruences in a covering $\mathcal{C}$ down to a minimal covering $\mathcal{C}_m$ in a finite number of steps, however this representation is not necessarily unique; that is, any covering system can be reduced to a minimal covering in at least one way.\\ 
\indent To our knowledge, the idea of reducing coverings was introduced formally by Krukenberg and honed by the first author and O. Trifonov in a $2022$ paper in which they proved the nonexistence of a distinct covering system with all moduli in the interval $[m, 8m]$ for $m\geq3$. Using their ideas, we extend their method for reducing coverings to a larger interval, namely $[m, 10m]$ for $m\geq3$.\\
\indent Here we list our first lemma (originally Corollary 7 in \cite{Dalton2}).
\begin{lem}\label{lem1}
    Let $\mathcal{C}$ be a covering such that $p^a|L$ for some prime $p$ and integer $a\geq1$. Suppose that there are $k$ congruences in $\mathcal{C}$ whose moduli are divisible by $p^a$. Then, if $k<p$, we can discard from $\mathcal{C}$ all congruences whose moduli are divisible by $p^a$ and still have a covering.
\end{lem}
Next, Krukenberg proved this lemma (originally a corollary),
\begin{lem}[Krukenberg]\label{lem2}
    Let $\mathcal{C}$ be a distinct covering with all moduli in the interval $[c, d]$. If $p$ is a prime and $a\geq1$ such that $p^a(p+1)>d$, then we can discard all congruences whose moduli are multiples of $p^a$ and still have a covering.
\end{lem}
Additionally, Krukenberg stated another lemma (originally a corollary) which we make use of,
\begin{lem}[Krukenberg]\label{lem3}
    Let $\mathcal{C}$ be a covering such that $p^a||L$ for some prime $p$ and integer $a\geq1$. Let $\mathcal{C}_1$ be a subset of $\mathcal{C}$ consisting of congruences whose moduli are divisible by $p^a$. Suppose $|\mathcal{C}_1|=p$ and the moduli of the congruences in $\mathcal{C}_1$ are $p^am_1,\ldots, p^am_p$. Then,\\
    (i) all the congruences in $\mathcal{C}_1$ can be replaced by a single congruence with modulus
    \begin{equation*}
        p^{a-1}\lcm(m_1,\ldots, m_p)
    \end{equation*}
    and the resulting set will still be a covering.\\
    (ii) if two of the above $p$ congruences are in the same class modulo $p^a$ we can discard all $p$ congruences and the resulting set will still be a covering.
\end{lem}
Lastly, we make use of one last lemma stated as follows,
\begin{lem}\label{lem4}
    Let $\mathcal{C}$ be a covering system and let $p$ be a prime. Let $\mathcal{C}_0$ be the subset of $\mathcal{C}$ whose moduli are not divisible by $p$. Let $\mathcal{M}_0$ be the list of moduli of the congruences in $\mathcal{C}_0$. Similarly, let $\mathcal{C}_1$ be the subset of $\mathcal{C}$ whose moduli are divisible by $p$, and let $\mathcal{M}_1$ be the list of the moduli of the congruences in $\mathcal{C}_1$.\\
    Reducing the covering $\mathcal{C}$ modulo $p$ produces $p$ coverings where\\
    \hspace{\parindent} (i) each modulus in $\mathcal{M}_0$ is used in each of the $p$ coverings but each modulus $n$ in $\mathcal{M}_1$ is replaced by $n/p$ and is used in just one of the $p$ coverings, and\\
    (ii) if two congruences in $\mathcal{C}_0$ are in the same class modulo a positive integer $q$, then after reducing they are in the same class modulo $q$ in each of the $p$ coverings; furthermore, if two congruences in $\mathcal{C}_0$ are not in the same class modulo a positive integer $q$, then after reducing they are not in the same class modulo $q$.
\end{lem}
\subsection*{Notation}
Throughout this paper, given $n\in\N$, we write $P(n)$ for its largest prime divisor with the convention that $P(1)=1$. Furthermore, adopting the notation of the first author and O. Trifonov, given a covering $\mathcal{C}$, we denote the least modulus $n_1$ by $m$, the largest modulus $n_k$ by $M$, and the least common multiple of all moduli in the covering by $L(\mathcal{C})=L$.

\section{Proof of Theorem \ref{thm2}}\label{sec:proofs}

All that is left for us is to prove Theorem \ref{thm2}. We note, for the exceptional $m$-values for which the sum of the reciprocals of the moduli is greater than $1$ and with $P(n)<\sqrt{9m+1}$, we break this proof up into three cases: reduction using Lemma \ref{lem3}, reduction using Lemma \ref{lem4}, and reduction using both lemmas. We now begin the proof.
\begin{proof}
    Assume for some integer $m\geq3$ there exists a distinct covering $\mathcal{C}$ with all moduli in the interval $[m,10m]$ and let $\mathcal{C}_m$ be a subset of $\mathcal{C}$. Consider the least common multiple $L$ of the moduli of the congruences in ${\mathcal{C}}_m$. By Lemma \ref{lem1}, if $p^a\vert L$ for some prime $p$ and a positive integer $a$, then the interval $[m,10m]$ contains at least $p$ multiples of $p^a$ that are not multiples of $p^{a+1}$. Since one of every $p$ consecutive multiples of $p^a$ is divisible by $p^{a+1}$, we deduce that the interval $[m,10m]$ contains at least $p+1$ multiples of $p^a$.\\
\indent Denote by $\mathcal{M}\subseteq[m,10m]$ the set of moduli from the congruences in ${\mathcal{C}}_m$. Let $p\geq\sqrt{9m+1}$ be a prime. The number of multiples of $p$ in the interval $[m,10m]$ is
\[
n_p:=\left\lfloor{\frac{10m}{p}}\right\rfloor-\left\lfloor{\frac{m-1}{p}}\right\rfloor=\frac{9m+1}{p}-\left\{\frac{10m}{p}\right\}+\left\{\frac{m-1}{p}\right\},
\]
where $\{x\}$ denotes the fractional part of $x$. Since for each $x$, $0\leq\{x\}<1$, we get
\[
n_p<\frac{9m+1}{p}+1\leq\sqrt{9m+1}+1\leq p+1.
\]
Thus, for each $p\geq\sqrt{9m+1}$, there are less than $p+1$ multiples of $p$ in the interval $[m,10m]$. Therefore, if $n$ is a modulus of one of the congruences in ${\mathcal{C}}_m$ (that is $n\in\mathcal{M}$), then all the prime divisors of $n$ are less than $\sqrt{9m+1}$. Since the density of integers covered by a congruence modulo $n$ is $1/n$ and ${\mathcal{C}}_m$ is a covering, we get
\begin{equation}
\sum_{\substack{m\leq n\leq 10m \\ P(n)<\sqrt{9m+1}}}\frac{1}{n} \geq \sum_{n\in \mathcal{M}} \frac{1}{n}\geq 1
\end{equation}
where $P(n)$ denotes the largest prime divisor of $n$.\\ \vspace{.2cm}
Let
\[
S_m=\sum_{n\in\mathcal{M}}\frac{1}{n}\ \ \ \text{and} \ \ \ T_m=\sum_{\substack{m\leq n\leq 10m \\ P(n)<\sqrt{9m+1}}}\frac{1}{n}.
\]
\indent We checked by direct computation that $T_m<1$ for all $m\in[117, 616000]$.\\
\indent Here we provide some details on how we showed that $T_m<1$ for all $m$ in $[117, 616000]$.\\
\indent First, via the Sieve of Eratosthenes we computed and stored $P(n)$, which denotes the largest prime divisor of $n$, for all $n$ from $2$ to $10\cdot616000$.\\
\indent Next note that $T_{m-1}\leq T_m+a_{m-1}$, where we define $a_{m-1}$ to be $\frac{1}{m-1}$ when $P(m-1)<\sqrt{9m-8}$, and we define $a_{m-1}$ to be $0$ when $P(m-1)\geq\sqrt{9m-8}$.\\ %this might look better as a piecewise function
\indent Indeed,
\begin{align*}
    T_{m-1}&=\sum_{\substack{m-1\leq n\leq 10m-10 \\ P(n)<\sqrt{9m-8}}}\frac{1}{n} \\
    &=a_{m-1}+\sum_{\substack{m\leq n\leq 10m-10 \\ P(n)<\sqrt{9m-8}}}\frac{1}{n} \\
    &\leq a_{m-1}+\sum_{\substack{m\leq n\leq 10m \\ P(n)<\sqrt{9m+1}}}\frac{1}{n} \\
    &=a_{m-1}+T_m.
\end{align*}
\indent So, we computed $T_{616000}$, and then using the inequality $T_{m-1}\leq T_m+a_{m-1}$ we get that $T_{615999}\leq T_{616000}+a_{615999}$. Iterating this method, we backtracked down to the last value of $m$ where the sum is less than $1$, which got us to $369082$. Thus, because
\[
    \sum_{m=369082}^{615999}a_m+T_{616000}=0.9999991593207759\ldots<1,
\]
we get that $T_m<1$ for all $m\in[369082, 616000]$. Next, we computed $T_{369081}$ and backtracked again; we got
\[
\sum_{m=224908}^{369081}a_m+T_{369082}=0.9999959672962349\ldots<1.
\]
Through these jumps, we confirmed that $T_m<1$ for each $m$ in the intervals,\\
\hfill \break
Table $3.1$ \ \ Intervals for $m\in[117, 616000]$ where $T_m<1$
\begin{center}
\begin{tabular}{|c|c|c|c|c|c|}
\hline
     $m$ & $[369082, 616000]$ & $[224908, 369081]$ & $[137109, 224907]$ & $[86560, 137108]$ & $[55427, 86559]$\\
     \hline
     $\ $ & $[35839, 55426]$ & $[23126, 35838]$ & $[15524, 23125]$ & $[10521, 15523]$ & $[7262, 10520]$\\
     \hline
     \hline
     $\ $ & $[5076, 7261]$ & $[3650, 5075]$ & $[2727, 3649]$ & $[2051, 2726]$ & $[1520, 2050]$\\
     \hline
     $\ $ & $[1225, 1519]$ & $[973, 1224]$ & $[761, 972]$ & $[647, 760]$ & $[537, 646]$\\
     \hline
     \hline
     $\ $ & $[436, 536]$ & $[373, 435]$ & $[302, 372]$ & $[276, 301]$ & $[247, 275]$\\
     \hline
     $\ $ & $[235, 246]$ & $[231, 234]$ & $[209, 230]$ & $[199, 208]$ & $[196, 198]$\\
     \hline
     \hline
     $\ $ & $[193, 195]$ & $[175, 192]$ & $[153, 174]$ & $[145, 152]$ & $[141, 144]$\\
     \hline
     $\ $ & $[126, 140]$ & $[120, 125]$ & $[118, 119]$ & $[117, 117]$ & $\ $\\
     \hline
\end{tabular}
\vspace*{-30pt}
\end{center}
\raggedright We had to compute $39$ values of $T_m$ to get to $m=117$. We then computed directly all values of $T_m$ for $m\in[3,116]$.
Since Balister et al. showed that the minimum modulus of a distinct covering system does not exceed $616000$, Theorem $\ref{thm2}$ holds when $m\geq117$.\\ 
\indent Also, by Theorem \ref{thm1}, since Krukenberg showed that there is no distinct covering system with moduli in $[3, 35]$, Theorem $\ref{thm2}$ holds when $m=3$.\\
\indent Furthermore, also by Theorem \ref{thm1}, there is no distinct covering system with moduli in $[4, 59]$; therefore, Theorem $\ref{thm2}$ holds when $m=4$ and $5$.\\
\indent There are $77$ occasions when $m\in[6, 116]$ and $T_m\geq1$. They are shown below.\\
\hfill \break
Table $3.2$ \ \ Values of $m\in[6, 116]$ where $T_m\geq1$
\begin{center}
\begin{tabular}{|c|c|c|c|c|c|}
\hline
     $m$ & $6$ & $7$ & $8$ & $9$ & $10$\\
     \hline
     $T_m$ & $1.5434041\ldots$ & $1.4225212\ldots$ & $1.3193863\ldots$ & $1.2297478\ldots$ & $1.1492575\ldots$\\
     \hline
     \hline
     $m$ & $11$ & $12$ & $13$ & $14$ & $15$\\
     \hline
     $T_m$ & $1.0680405\ldots$ & $1.0853025\ldots$ & $1.0257181\ldots$ & $1.2314694\ldots$ & $1.1804546\ldots$\\
     \hline
     \hline
     $m$ & $16$ & $17$ & $18$ & $19$ & $20$\\
     \hline
     $T_m$ & $1.1265315\ldots$ & $1.0822173\ldots$ & $1.0991689\ldots$ & $1.2221015\ldots$ & $1.2475906\ldots$\\
     \hline
     \hline
     $m$ & $21$ & $22$ & $23$ & $24$ & $25$\\
     \hline
     $T_m$ & $1.2071602\ldots$ & $1.1687162\ldots$ & $1.1321704\ldots$ & $1.1449396\ldots$ & $1.1196020\ldots$\\
     \hline
     \hline
     $m$ & $26$ & $27$ & $28$ & $29$ & $30$\\
     \hline
     $T_m$ & $1.0913227\ldots$ & $1.0603527\ldots$ & $1.0341865\ldots$ & $1.0054409\ldots$ & $1.0155426\ldots$\\
     \hline
     \hline
     $m$ & $33$ & $34$ & $35$ & $36$ & $37$\\
     \hline
     $T_m$ & $1.1196755\ldots$ & $1.0982484\ldots$ & $1.0746092\ldots$ & $1.0573066\ldots$ & $1.0350309\ldots$\\
     \hline
     \hline
     $m$ & $38$ & $39$ & $40$ & $41$ & $42$\\
     \hline
     $T_m$ & $1.0430168\ldots$ & $1.0507825\ldots$ & $1.0327178\ldots$ & $1.1255515\ldots$ & $1.1327286\ldots$\\
     \hline
     \hline
     $m$ & $43$ & $44$ & $45$ & $46$ & $47$\\
     \hline
     $T_m$ & $1.1136031\ldots$ & $1.1181906\ldots$ & $1.1044477\ldots$ & $1.0887949\ldots$ & $1.0930962\ldots$\\
     \hline
     \hline
     $m$ & $48$ & $49$ & $50$ & $51$ & $52$\\
     \hline
     $T_m$ & $1.0993856\ldots$ & $1.0847168\ldots$ & $1.0703532\ldots$ & $1.0562705\ldots$ & $1.0424881\ldots$\\
     \hline
     \hline
     $m$ & $53$ & $54$ & $55$ & $56$ & $59$\\
     \hline
     $T_m$ & $1.0270561\ldots$ & $1.0326429\ldots$ & $1.0196123\ldots$ & $1.0032162\ldots$ & $1.0824342\ldots$\\
     \hline
     \hline
     $m$ & $60$ & $61$ & $62$ & $63$ & $64$\\
     \hline
     $T_m$ & $1.0891372\ldots$ & $1.0757682\ldots$ & $1.0790256\ldots$ & $1.0870206\ldots$ & $1.0742800\ldots$\\
     \hline
     \hline
     $m$ & $65$ & $66$ & $67$ & $68$ & $69$\\
     \hline
     $T_m$ & $1.0648374\ldots$ & $1.0509680\ldots$ & $1.0388285\ldots$ & $1.0447480\ldots$ & $1.0344111\ldots$\\
     \hline
     \hline
     $m$ & $70$ & $71$ & $72$ & $73$ & $74$\\
     \hline
     $T_m$ & $1.0227899\ldots$ & $1.0113491\ldots$ & $1.0155372\ldots$ & $1.0071561\ldots$ & $1.0098754\ldots$\\
     \hline
     \hline
     $m$ & $75$ & $76$ & $95$ & $107$ & $108$\\
     \hline
     $T_m$ & $1.0138951\ldots$ & $1.0045179\ldots$ & $1.0010795\ldots$ & $1.0253732\ldots$ & $1.0281605\ldots$\\
     \hline
     \hline
     $m$ & $109$ & $110$ & $111$ & $112$ & $113$\\
     \hline
     $T_m$ & $1.0225837\ldots$ & $1.0244085\ldots$ & $1.0180358\ldots$ & $1.0198247\ldots$ & $1.0135636\ldots$\\
     \hline
     \hline
     $m$ & $114$ & $115$ & $116$ & &\\
     \hline
     $T_m$ & $1.0162068\ldots$ & $1.0091786\ldots$ & $1.0039439\ldots$ & &\\
     \hline
\end{tabular}
\end{center}
\indent So far, we have used Lemma \ref{lem1} only with $a=1$. Next, we use Lemma \ref{lem2} for all $a\geq1$.\\
\hfill \break
\indent Define
\vspace{-20pt}
\begin{center}
    \[ L_m = \begin{cases} 
    5040=2^4\cdot3^2\cdot5\cdot7\ \text{if}\ m\in\{6, 7, 8, 9\}\\
    
    10080=2^5\cdot3^2\cdot5\cdot7\ \text{if}\ m=10\\
    
    30240=2^5\cdot3^3\cdot5\cdot7\ \text{if}\ m\in\{11, 12, 13, 14\}\\
    
    1663200=2^5\cdot3^3\cdot5^2\cdot7\cdot11\ \text{if}\ m\in[15, 19]\\
    43243200=2^6\cdot3^3\cdot5^2\cdot7\cdot11\cdot13\ \text{if}\ m\in[20, 30]\\
    2.2054032\times10^9=2^6\cdot3^4\cdot5^2\cdot7\cdot11\cdot13\cdot17\ \text{if}\ m\in[33, 38]\\
    4.4108064\times10^9=2^7\cdot3^4\cdot5^2\cdot7\cdot11\cdot13\cdot17\ \text{if}\ m=39\\
    3.08756448\times10^{10}=2^7\cdot3^4\cdot5^2\cdot7^2\cdot11\cdot13\cdot17\ \text{if}\ m\in\{40, 41\}\\
    5.866372512\times10^{11}=2^7\cdot3^4\cdot5^2\cdot7^2\cdot11\cdot13\cdot17\cdot19\ \text{if}\ m\in[42, 56]\cup\{59\}\\
    1.3492656778\times10^{13}=2^7\cdot3^4\cdot5^2\cdot7^2\cdot11\cdot13\cdot17\cdot19\cdot23\ \text{if}\ m\in[60, 74]\\
    6.7463283888\times10^{13}=2^7\cdot3^4\cdot5^3\cdot7^2\cdot11\cdot13\cdot17\cdot19\cdot23\ \text{if}\ m\in\{75, 76\}\\
    1.3492656778\times10^{14}=2^8\cdot3^4\cdot5^3\cdot7^2\cdot11\cdot13\cdot17\cdot19\cdot23\ \text{if}\ m=95\\
    1.1738611937\times10^{16}=2^8\cdot3^5\cdot5^3\cdot7^2\cdot11\cdot13\cdot17\cdot19\cdot23\cdot29\ \text{if}\ m\in\{107, 108\}\\
    3.6389695329\times10^{17}=2^8\cdot3^5\cdot5^3\cdot7^2\cdot11\cdot13\cdot17\cdot19\cdot23\cdot29\cdot31\ \text{if}\ m\in[109, 116]\\
       \end{cases}
    \]
\end{center}
By Lemma 2.2 one checks directly that $L$ divides $L_m$ for each $m$ which it is defined above. Then for such $m\notin\{6, 7, 8, 9, 15, 16, 18, 20, 21, 22, 23, 24, 25, 33, 42, 43, 44, 45, 46, 47, 48, 49, 50\}$, a direct computation shows we have
\begin{equation}
\sum_{n\in \mathcal{M}} \frac{1}{n}\leq\sum_{\substack{d|L_m \\ m\leq d\leq10m}}\frac{1}{d}<1.
\end{equation}
\indent Since $(3.1)$ contradicts the second inequality in 
$(3.2)$, the proof is complete for all but $m\in\{6, 7, 8, 9, 15, 16, 18, 20, 21, 22, 23, 24, 25, 33, 42, 43, 44, 45, 46, 47, 48, 49, 50\}$. We now split our proof into three cases: reduction using Lemma \ref{lem3}, reduction using Lemma \ref{lem4}, and reduction using both lemmas.
\subsection*{Reduction via Lemma \ref{lem3}}
For all $m\in\{6, 16, 20, 33, 42, 43, 45\}$, we can get the sums of the reciprocals of the moduli to be $<1$ by only using Lemma \ref{lem3}. We write these cases out here.\\
\hfill \break
\indent For $m=6$, suppose we have a distinct covering with the moduli of all the congruences in $[6, 60]$. We have,
\begin{equation*}
    \sum_{\substack{d|L_6 \\ 6\leq d\leq60}}\frac{1}{d}\approx1.3761905.
\end{equation*}
Since there are $7$ permissible multiples of $7$ in the interval $[6, 60]$, namely $7, 2\cdot7, 3\cdot7, 2^2\cdot7, 5\cdot7, 2\cdot3\cdot7$, and $2^3\cdot7$, we can replace the congruences with these moduli by a single congruence modulo $120$ and still have a covering. This brings the sum of the reciprocals down to approximately $1.016667$.\\
\indent Further, since there are $2$ permissible multiples of $2^4$ in the interval $[6, 60]$, namely $2^4$ and $2^4\cdot3$, we can replace the congruences with these moduli by a single congruence modulo $2^3\cdot3$ and still have a covering. This brings the sum of the reciprocals down to approximately $0.933333<1$, and thus we have a contradiction. Therefore, the theorem holds for $m=6$.\\
\hfill \break
\indent For $m=16$, suppose we have a distinct covering with the moduli of all the congruences in $[16, 160]$. We have,
\begin{equation*}
    \sum_{\substack{d|L_{16} \\ 16\leq d\leq160}}\frac{1}{d}\approx1.0370689.
\end{equation*}
Since there are $5$ permissible multiples of $5^2$ in the interval $[16, 160]$, namely $5^2, 2\cdot5^2, 3\cdot5^2, 4\cdot5^2$, and $6\cdot5^2$, we can replace the congruences with these moduli by a single congruence modulo $5\cdot12$ and still have a covering. This brings the sum of the reciprocals down to approximately $1.030402$.\\
\indent Further, since there are $11$ permissible multiples of $11$ in the interval $[16, 160]$, namely $2\cdot11, 3\cdot11, 4\cdot11, 5\cdot11, 6\cdot11, 7\cdot11, 8\cdot11, 9\cdot11, 10\cdot11, 12\cdot11$ and $14\cdot11$, we can replace the congruences with these moduli by a single congruence modulo $2520$ and still have a covering. This brings the sum of the reciprocals down to approximately $0.841369<1$, and thus we have a contradiction. Therefore, the theorem holds for $m=16$.\\
\hfill \break
For $m=20$, suppose we have a distinct covering with the moduli of all the congruences in $[16, 160]$. We have,
\begin{equation*}
    \sum_{\substack{d|L_{20} \\ 20\leq d\leq200}}\frac{1}{d}\approx1.1565610.
\end{equation*}
Since there are $13$ permissible multiples of $13$ in the interval $[20, 200]$, namely $2\cdot13, 3\cdot13,\ldots, 12\cdot13, 14\cdot13$, and $15\cdot13$, we can replace the congruences with these moduli by a single congruence modulo $27720$ and still have a covering. This brings the sum of the reciprocals down to approximately $0.984189<1$, and thus we have a contradiction. Therefore, the theorem holds for $m=20$.\\
\hfill \break
For $m=33$, suppose we have a distinct covering with the moduli of all the congruences in $[33, 330]$. We have,
\begin{equation*}
    \sum_{\substack{d|L_{33} \\ 33\leq d\leq330}}\frac{1}{d}\approx 1.0200675.
\end{equation*}
Since there are $17$ permissible multiples of $17$ in the interval $[33, 330]$, namely $2\cdot17, 3\cdot17,\ldots, 16\cdot17$, and $18\cdot13$, we can replace the congruences with these moduli by a single congruence modulo $720720$ and still have a covering. This brings the sum of the reciprocals down to approximately $0.876758<1$, and thus we have a contradiction. Therefore, the theorem holds for $m=33$.\\
\hfill \break
\indent For $m=42$, suppose we have a distinct covering with the moduli of all the congruences in $[42, 420]$. We have,
\begin{equation*}
    \sum_{\substack{d|L_{42} \\ 42\leq d\leq420}}\frac{1}{d}\approx1.0768676.
\end{equation*}
\indent Since there are $19$ permissible multiples of $19$ in the interval $[42, 420]$, namely $19\cdot3, 19\cdot4,\ldots, 19\cdot18, 19\cdot20, 19\cdot21,$ and $19\cdot22$, we can replace the congruences with these moduli by a single congruence modulo $12252240$ and still have a covering. This brings the sum of the reciprocals down to approximately $0.964332<1$, and thus we have a contradiction. Therefore, the theorem holds for $m=42$.\\
\hfill \break
For $m=43$, suppose we have a distinct covering with the moduli of all the congruences in $[43, 430]$. We have,
\begin{equation*}
    \sum_{\substack{d|L_{43} \\ 43\leq d\leq430}}\frac{1}{d}\approx 1.0577420.
\end{equation*}
Since there are $19$ permissible multiples of $19$ in the interval $[43, 430]$, namely $3\cdot19, 4\cdot19,\ldots, 18\cdot19, 20\cdot19, 21\cdot19$, and $22\cdot19$, we can replace the congruences with these moduli by a single congruence modulo $12252240$ and still have a covering. This brings the sum of the reciprocals down to approximately $0.945206<1$, and thus we have a contradiction. Therefore, the theorem holds for $m=43$.\\
\hfill \break
For $m=45$, suppose we have a distinct covering with the moduli of all the congruences in $[45, 450]$. We have,
\begin{equation*}
    \sum_{\substack{d|L_{45} \\ 45\leq d\leq450}}\frac{1}{d}\approx 1.0485866.
\end{equation*}
Since there are $19$ permissible multiples of $19$ in the interval $[45, 450]$, namely $3\cdot19, 4\cdot19,\ldots, 18\cdot19, 20\cdot19, 21\cdot19$, and $22\cdot19$, we can replace the congruences with these moduli by a single congruence modulo $12252240$ and still have a covering. This brings the sum of the reciprocals down to approximately $0.936051<1$, and thus we have a contradiction. Therefore, the theorem holds for $m=45$.
\subsection*{Reduction via Lemma \ref{lem4}}
For all $m\in\{8, 9, 15, 18, 22, 23, 24, 25, 44, 46, 47, 48, 49, 50\}$, we can get the sums of the reciprocals of the moduli to be $<1$ by only using Lemma \ref{lem4}. We write these cases out here.\\
\hfill \break
\indent For $m=8$, suppose we have a distinct covering with the moduli of all the congruences in $[8, 80]$. We have,
\begin{equation*}
    \sum_{\substack{d|L_8 \\ 8\leq d\leq80}}\frac{1}{d}\approx1.1232142.
\end{equation*}
\indent Denote the resulting covering by $\mathcal{C}'$ and denote the list of the moduli of the congruences in $\mathcal{C}'$ by $\mathcal{M}'$.\\
\indent Next, we reduce the covering $\mathcal{C}'$ modulo $7$ and obtain $7$ coverings, say $\mathcal{C}'_0, \mathcal{C}'_1,\ldots, \mathcal{C}'_6$, whose moduli are $\mathcal{M}'_0, \mathcal{M}'_1,\ldots, \mathcal{M}'_6$, respectively. The moduli in
\begin{equation*}
    \mathcal{M}_0=\{8, 9, 10, 12, 15, 16, 18, 20, 24, 30, 36, 40, 45, 48, 60, 72, 80\}
\end{equation*}
can be used in all $7$ coverings, and each modulus in
\begin{equation*}
    \mathcal{M}_1=\{2, 3, 4, 5, 6, 8, 9, 10\}
\end{equation*}
can be used in just one of the coverings.\\
\indent Now,
\begin{equation*}
    \sum_{d\in\mathcal{M}_0}\frac{1}{d}\approx0.868056<1,
\end{equation*}
which means that each of the seven coverings needs at least one of the $8$ moduli in $\mathcal{M}_1$ to get the sum of the reciprocals to be $1$ or more.\\
\indent After relabeling, we can assume that $2\in\mathcal{M}'_0,\ 3\in\mathcal{M}'_1,\ldots,\ \text{and}\ 6\in\mathcal{M}'_4$. This leaves us with three moduli, $8, 9, 10$, to distribute into the remaining two coverings, but none of these remaining moduli alone are enough to get the sum of the reciprocals to be $1$ or more, and thus we have a contradiction. Therefore, the theorem holds for $m=8$.\\
\hfill \break
\indent For $m=9$, suppose we have a distinct covering with the moduli of all the congruences in $[9, 90]$. We have,
\begin{equation*}
    \sum_{\substack{d|L_9 \\ 9\leq d\leq90}}\frac{1}{d}\approx1.0212301.
\end{equation*}
\indent Denote the resulting covering by $\mathcal{C}'$ and denote the list of the moduli of the congruences in $\mathcal{C}'$ by $\mathcal{M}'$.\\
\indent Next, we reduce the covering $\mathcal{C}'$ modulo $7$ and obtain $7$ coverings, say $\mathcal{C}'_0, \mathcal{C}'_1,\ldots, \mathcal{C}'_6$, whose moduli are $\mathcal{M}'_0, \mathcal{M}'_1,\ldots, \mathcal{M}'_6$, respectively. The moduli in
\begin{equation*}
    \mathcal{M}_0=\{9, 10, 12, 15, 16, 18, 20, 24, 30, 36, 40, 45, 48, 60, 72, 80, 90\}
\end{equation*}
can be used in all $7$ coverings, and each modulus in
\begin{equation*}
    \mathcal{M}_1=\{2, 3, 4, 5, 6, 8, 9, 10, 12\}
\end{equation*}
can be used in just one of the coverings.\\
\indent Now,
\begin{equation*}
    \sum_{d\in\mathcal{M}_0}\frac{1}{d}\approx0.754167<1,
\end{equation*}
which means that each of the seven coverings needs at least one of the $9$ moduli in $\mathcal{M}_1$ to get the sum of the reciprocals to be $1$ or more.\\
\indent After relabeling, we can assume that $2\in\mathcal{M}'_0,\ 3\in\mathcal{M}'_1,\ \text{and}\ 4\in\mathcal{M}'_2$. This leaves us with six moduli to distribute into the remaining four coverings, but none of these remaining moduli alone are enough to get the sum of the reciprocals to be $1$ or more, and thus we have a contradiction. Therefore, the theorem holds for $m=9$.\\
\hfill \break
\indent For $m=15$, suppose we have a distinct covering with the moduli of all the congruences in $[15, 150]$. We have,
\begin{equation*}
    \sum_{\substack{d|L_{15} \\ 15\leq d\leq150}}\frac{1}{d}\approx1.0909921.
\end{equation*}
\indent Denote the resulting covering by $\mathcal{C}'$ and denote the list of the moduli of the congruences in $\mathcal{C}'$ by $\mathcal{M}'$.\\
\indent Next, we reduce the covering $\mathcal{C}'$ modulo $11$ and obtain $11$ coverings, say $\mathcal{C}'_0, \mathcal{C}'_1,\ldots, \mathcal{C}'_{10}$, whose moduli are $\mathcal{M}'_0, \mathcal{M}'_1,\ldots, \mathcal{M}'_{10}$, respectively. The moduli in
\[
    \mathcal{M}_0=\{15, 16,\ldots, 150\}
\]
can be used in all $11$ coverings, and each modulus in
\begin{equation*}
    \mathcal{M}_1=\{2, 3, 4, 5, 6, 7, 8, 9, 10, 12\}
\end{equation*}
can be used in just one of the coverings.\\
\indent Now,
\begin{equation*}
    \sum_{d\in\mathcal{M}_0}\frac{1}{d}\approx0.908056<1,
\end{equation*}
which means that each of the eleven coverings needs at least one of the $10$ moduli in $\mathcal{M}_1$ to get the sum of the reciprocals to be $1$ or more.\\
\indent After relabeling, we can assume that $2\in\mathcal{M}'_0,\ 3\in\mathcal{M}'_1,\ldots, \ \text{and}\ 12\in\mathcal{M}'_9$. This leaves us with no moduli to distribute into the last covering, and therefore the sum of the reciprocals cannot be $1$ or more, and thus we have a contradiction. Therefore, the theorem holds for $m=15$.\\
\hfill \break
\indent For $m=18$, suppose we have a distinct covering with the moduli of all the congruences in $[18, 180]$. We have,
\begin{equation*}
    \sum_{\substack{d|L_{18} \\ 18\leq d\leq180}}\frac{1}{d}\approx1.0035335.
\end{equation*}
\indent Denote the resulting covering by $\mathcal{C}'$ and denote the list of the moduli of the congruences in $\mathcal{C}'$ by $\mathcal{M}'$.\\
\indent Next, we reduce the covering $\mathcal{C}'$ modulo $11$ and obtain $11$ coverings, say $\mathcal{C}'_0, \mathcal{C}'_1,\ldots, \mathcal{C}'_{10}$, whose moduli are $\mathcal{M}'_0, \mathcal{M}'_1,\ldots, \mathcal{M}'_{10}$, respectively. The moduli in
\[
    \mathcal{M}_0=\{18, 20,\ldots, 180\}
\]
can be used in all $11$ coverings, and each modulus in
\begin{equation*}
    \mathcal{M}_1=\{2, 3, 4, 5, 6, 7, 8, 9, 10, 12, 14, 15, 16\}
\end{equation*}
can be used in just one of the coverings.\\
\indent Now,
\begin{equation*}
    \sum_{d\in\mathcal{M}_0}\frac{1}{d}\approx0.802361<1,
\end{equation*}
which means that each of the eleven coverings needs at least one of the $13$ moduli in $\mathcal{M}_1$ to get the sum of the reciprocals to be $1$ or more.\\
\indent After relabeling, we can assume that $2\in\mathcal{M}'_0,\ 3\in\mathcal{M}'_1,\ldots, \ \text{and}\ 5\in\mathcal{M}'_3$. This leaves us with eight moduli, $6, 7, 8, 9, 10, 12, 14, 15, 16$, to distribute into the remaining seven coverings, but none of these remaining moduli alone are enough to get the sum of the reciprocals to be $1$ or more, and thus we have a contradiction. Therefore, the theorem holds for $m=18$.\\
\hfill \break
For $m=22$, suppose we have a distinct covering with the moduli of all the congruences in $[22, 220]$. We have,
\begin{equation*}
    \sum_{\substack{d|L_{22} \\ 22\leq d\leq220}}\frac{1}{d}\approx1.0776866.
\end{equation*}
\indent Denote the resulting covering by $\mathcal{C}'$ and denote the list of the moduli of the congruences in $\mathcal{C}'$ by $\mathcal{M}'$.\\
\indent Next, we reduce the covering $\mathcal{C}'$ modulo $13$ and obtain $13$ coverings, say $\mathcal{C}'_0, \mathcal{C}'_1,\ldots, \mathcal{C}'_{12}$, whose moduli are $\mathcal{M}'_0, \mathcal{M}'_1,\ldots, \mathcal{M}'_{12}$, respectively. The moduli in
\[
    \mathcal{M}_0=\{22, 24,\ldots, 220\}
\]
can be used in all $13$ coverings, and each modulus in
\begin{equation*}
    \mathcal{M}_1=\{2, 3, 4, 5, 6, 7, 8, 9, 10, 11, 12, 14, 15, 16\}
\end{equation*}
can be used in just one of the coverings.\\
\indent Now,
\begin{equation*}
    \sum_{d\in\mathcal{M}_0}\frac{1}{d}\approx0.900471<1,
\end{equation*}
which means that each of the thirteen coverings needs at least one of the $14$ moduli in $\mathcal{M}_1$ to get the sum of the reciprocals to be $1$ or more.\\
\indent After relabeling, we can assume that $2\in\mathcal{M}'_0,\ 3\in\mathcal{M}'_1,\ldots,\ \text{and}\ 10\in\mathcal{M}'_{8}$. This leaves us with five moduli, $11, 12, 14, 15, 16$, to distribute into the remaining four coverings, but none of these remaining moduli alone are enough to get the sum of the reciprocals to be $1$ or more, and thus we have a contradiction. Therefore, the theorem holds for $m=22$.\\
\hfill \break
\indent For $m=23$, suppose we have a distinct covering with the moduli of all the congruences in $[23, 230]$. We have,
\begin{equation*}
    \sum_{\substack{d|L_{23} \\ 23\leq d\leq230}}\frac{1}{d}\approx1.0411408.
\end{equation*}
\indent Denote the resulting covering by $\mathcal{C}'$ and denote the list of the moduli of the congruences in $\mathcal{C}'$ by $\mathcal{M}'$.\\
\indent Next, we reduce the covering $\mathcal{C}'$ modulo $13$ and obtain $13$ coverings, say $\mathcal{C}'_0, \mathcal{C}'_1,\ldots, \mathcal{C}'_{13}$, whose moduli are $\mathcal{M}'_0, \mathcal{M}'_1,\ldots, \mathcal{M}'_{12}$, respectively. The moduli in
\[
    \mathcal{M}_0=\{24, 25,\ldots, 225\}
\]
can be used in all $13$ coverings, and each modulus in
\begin{equation*}
    \mathcal{M}_1=\{2, 3, 4, 5, 6, 7, 8, 9, 10, 11, 12, 14, 15, 16\}
\end{equation*}
can be used in just one of the coverings.\\
\indent Now,
\begin{equation*}
    \sum_{d\in\mathcal{M}_0}\frac{1}{d}\approx0.863925<1,
\end{equation*}
which means that each of the thirteen coverings needs at least one of the $14$ moduli in $\mathcal{M}_1$ to get the sum of the reciprocals to be $1$ or more.\\
\indent After relabeling, we can assume that $2\in\mathcal{M}'_0,\ 3\in\mathcal{M}'_1,\ldots, \ \text{and}\ 7\in\mathcal{M}'_5$. This leaves us with eight moduli, $8, 9, 10, 11, 12, 14, 15, 16$, to distribute into the remaining seven coverings, but none of these remaining moduli alone are enough to get the sum of the reciprocals to be $1$ or more, and thus we have a contradiction. Therefore, the theorem holds for $m=23$.\\
\hfill \break
\indent For $m=24$, suppose we have a distinct covering with the moduli of all the congruences in $[24, 240]$. We have,
\begin{equation*}
    \sum_{\substack{d|L_{24} \\ 24\leq d\leq240}}\frac{1}{d}\approx1.0539099.
\end{equation*}
\indent Denote the resulting covering by $\mathcal{C}'$ and denote the list of the moduli of the congruences in $\mathcal{C}'$ by $\mathcal{M}'$.\\
\indent Next, we reduce the covering $\mathcal{C}'$ modulo $13$ and obtain $13$ coverings, say $\mathcal{C}'_0, \mathcal{C}'_1,\ldots, \mathcal{C}'_{13}$, whose moduli are $\mathcal{M}'_0, \mathcal{M}'_1,\ldots, \mathcal{M}'_{12}$, respectively. The moduli in
\[
    \mathcal{M}_0=\{24, 25,\ldots, 240\}
\]
can be used in all $13$ coverings, and each modulus in
\begin{equation*}
    \mathcal{M}_1=\{2, 3, 4, 5, 6, 7, 8, 9, 10, 11, 12, 14, 15, 16, 18\}
\end{equation*}
can be used in just one of the coverings.\\
\indent Now,
\begin{equation*}
    \sum_{d\in\mathcal{M}_0}\frac{1}{d}\approx0.872421<1,
\end{equation*}
which means that each of the thirteen coverings needs at least one of the $15$ moduli in $\mathcal{M}_1$ to get the sum of the reciprocals to be $1$ or more.\\
\indent After relabeling, we can assume that $2\in\mathcal{M}'_0,\ 3\in\mathcal{M}'_1,\ldots, \ \text{and}\ 7\in\mathcal{M}'_5$. This leaves us with nine moduli, $8, 9, 10, 11, 12, 14, 15, 16, 18$, to distribute into the remaining seven coverings, but none of these remaining moduli alone are enough to get the sum of the reciprocals to be $1$ or more, and thus we have a contradiction. Therefore, the theorem holds for $m=24$.\\
\hfill \break
\indent For $m=25$, suppose we have a distinct covering with the moduli of all the congruences in $[25, 250]$. We have,
\begin{equation*}
    \sum_{\substack{d|L_{25} \\ 25\leq d\leq250}}\frac{1}{d}\approx1.0122433.
\end{equation*}
\indent Denote the resulting covering by $\mathcal{C}'$ and denote the list of the moduli of the congruences in $\mathcal{C}'$ by $\mathcal{M}'$.\\
\indent Next, we reduce the covering $\mathcal{C}'$ modulo $13$ and obtain $13$ coverings, say $\mathcal{C}'_0, \mathcal{C}'_1,\ldots, \mathcal{C}'_{12}$, whose moduli are $\mathcal{M}'_0, \mathcal{M}'_1,\ldots, \mathcal{M}'_{12}$, respectively. The moduli in
\[
    \mathcal{M}_0=\{25, 27,\ldots, 240\}
\]
can be used in all $13$ coverings, and each modulus in
\begin{equation*}
    \mathcal{M}_1=\{2, 3, 4, 5, 6, 7, 8, 9, 10, 11, 12, 14, 15, 16, 18\}
\end{equation*}
can be used in just one of the coverings.\\
\indent Now,
\begin{equation*}
    \sum_{d\in\mathcal{M}_0}\frac{1}{d}\approx0.830754<1,
\end{equation*}
which means that each of the thirteen coverings needs at least one of the $15$ moduli in $\mathcal{M}_1$ to get the sum of the reciprocals to be $1$ or more.\\
\indent After relabeling, we can assume that $2\in\mathcal{M}'_0,\ 3\in\mathcal{M}'_1,\ldots, \ \text{and}\ 5\in\mathcal{M}'_3$. This leaves us with eleven moduli, $6, 7, 8, 9, 10, 11, 12, 14, 15, 16, 18$, to distribute into the remaining nine coverings, but none of these remaining moduli alone are enough to get the sum of the reciprocals to be $1$ or more, and thus we have a contradiction. Therefore, the theorem holds for $m=25$.\\
\hfill \break
\indent For $m=44$, suppose we have a distinct covering with the moduli of all the congruences in $[44, 440]$. We have,
\begin{equation*}
    \sum_{\substack{d|L_{44} \\ 44\leq d\leq440}}\frac{1}{d}\approx1.0623295.
\end{equation*}
\indent Next, denote the resulting covering by $\mathcal{C}'$ and denote the list of the moduli of the congruences in $\mathcal{C}'$ by $\mathcal{M}'$.\\
\indent Next, we reduce the covering $\mathcal{C}'$ modulo $19$ and obtain $19$ coverings, say $\mathcal{C}'_0, \mathcal{C}'_1,\ldots, \mathcal{C}'_{18}$, whose moduli are $\mathcal{M}'_0, \mathcal{M}'_1,\ldots, \mathcal{M}'_{18}$, respectively.
The moduli in
\begin{equation*}
    \mathcal{M}_0=\{44, 45,\ldots, 440\}
\end{equation*}
can be used in all $19$ coverings, and each modulus in
\begin{equation*}
    \mathcal{M}_1=\{3, 4,\ldots,18, 20,\ldots ,22\}
\end{equation*}
can be used in just one of the coverings.\\
\indent Now,
\begin{equation*}
    \sum_{d\in\mathcal{M}_0}\frac{1}{d}\approx0.949794<1,
\end{equation*}
which means that each of the nineteen coverings needs at least one of the $19$ moduli in $\mathcal{M}_1$ to get the sum of the reciprocals to be $1$ or more.\\
\indent After relabeling, we can assume that $3\in\mathcal{M}'_0,\ 4\in\mathcal{M}'_1,\ldots,\ \text{and}\ 18\in\mathcal{M}'_{15}$. This leaves us with three moduli, $20, 21, 22$, to distribute into the remaining three coverings, but none of these remaining moduli alone are enough to get the sum of the reciprocals to be $1$ or more, and thus we have a contradiction. Therefore, the theorem holds for $m=44$.\\
\hfill \break
\indent For $m=46$, suppose we have a distinct covering with the moduli of all the congruences in $[46, 460]$. We have,
\begin{equation*}
    \sum_{\substack{d|L_{46} \\ 46\leq d\leq460}}\frac{1}{d}\approx1.0329338.
\end{equation*}
\indent Next, denote the resulting covering by $\mathcal{C}'$ and denote the list of the moduli of the congruences in $\mathcal{C}'$ by $\mathcal{M}'$.\\
\indent Next, we reduce the covering $\mathcal{C}'$ modulo $19$ and obtain $19$ coverings, say $\mathcal{C}'_0, \mathcal{C}'_1,\ldots, \mathcal{C}'_{18}$, whose moduli are $\mathcal{M}'_0, \mathcal{M}'_1,\ldots, \mathcal{M}'_{18}$, respectively.
The moduli in
\begin{equation*}
    \mathcal{M}_0=\{48, 49,\ldots, 459\}
\end{equation*}
can be used in all $19$ coverings, and each modulus in
\begin{equation*}
    \mathcal{M}_1=\{3, 4,\ldots, 18, 20,\ldots, 22, 24\}
\end{equation*}
can be used in just one of the coverings.\\
\indent Now,
\begin{equation*}
    \sum_{d\in\mathcal{M}_0}\frac{1}{d}\approx0.918205<1,
\end{equation*}
which means that each of the nineteen coverings needs at least one of the $20$ moduli in $\mathcal{M}_1$ to get the sum of the reciprocals to be $1$ or more.\\
\indent After relabeling, we can assume that $3\in\mathcal{M}'_0,\ 4\in\mathcal{M}'_1,\ldots,\ \text{and}\ 12\in\mathcal{M}'_{9}$. This leaves us with ten moduli, $13, 14, 15, 16, 17, 18, 20, 21, 22, 24$, to distribute into the remaining nine coverings, but none of these remaining moduli alone are enough to get the sum of the reciprocals to be $1$ or more, and thus we have a contradiction. Therefore, the theorem holds for $m=46$.\\
\hfill \break
\indent For $m=47$, suppose we have a distinct covering with the moduli of all the congruences in $[47, 470]$. We have,
\begin{equation*}
    \sum_{\substack{d|L_{47} \\ 47\leq d\leq470}}\frac{1}{d}\approx1.0372351.
\end{equation*}
\indent Next, denote the resulting covering by $\mathcal{C}'$ and denote the list of the moduli of the congruences in $\mathcal{C}'$ by $\mathcal{M}'$.\\
\indent Next, we reduce the covering $\mathcal{C}'$ modulo $19$ and obtain $19$ coverings, say $\mathcal{C}'_0, \mathcal{C}'_1,\ldots, \mathcal{C}'_{18}$, whose moduli are $\mathcal{M}'_0, \mathcal{M}'_1,\ldots, \mathcal{M}'_{18}$, respectively.
The moduli in
\begin{equation*}
    \mathcal{M}_0=\{48, 49,\ldots, 468\}
\end{equation*}
can be used in all $19$ coverings, and each modulus in
\begin{equation*}
    \mathcal{M}_1=\{3, 4,\ldots,18, 20,\ldots, 24\}
\end{equation*}
can be used in just one of the coverings.\\
\indent Now,
\begin{equation*}
    \sum_{d\in\mathcal{M}_0}\frac{1}{d}\approx0.922506<1,
\end{equation*}
which means that each of the nineteen coverings needs at least one of the $22$ moduli in $\mathcal{M}_1$ to get the sum of the reciprocals to be $1$ or more.\\
\indent After relabeling, we can assume that $3\in\mathcal{M}'_0,\ 4\in\mathcal{M}'_1,\ldots,\ \text{and}\ 12\in\mathcal{M}'_9$. This leaves us with ten moduli, $13,\ldots, 18, 20,\ldots, 24$, to distribute into the remaining nine coverings, but none of these remaining moduli alone are enough to get the sum of the reciprocals to be $1$ or more, and thus we have a contradiction. Therefore, the theorem holds for $m=47$.\\
\hfill \break
\indent For $m=48$, suppose we have a distinct covering with the moduli of all the congruences in $[48, 480]$. We have,
\begin{equation*}
    \sum_{\substack{d|L_{48} \\ 48\leq d\leq480}}\frac{1}{d}\approx1.0435245.
\end{equation*}
\indent Next, denote the resulting covering by $\mathcal{C}'$ and denote the list of the moduli of the congruences in $\mathcal{C}'$ by $\mathcal{M}'$.\\
\indent Next, we reduce the covering $\mathcal{C}'$ modulo $19$ and obtain $19$ coverings, say $\mathcal{C}'_0, \mathcal{C}'_1,\ldots, \mathcal{C}'_{18}$, whose moduli are $\mathcal{M}'_0, \mathcal{M}'_1,\ldots, \mathcal{M}'_{18}$, respectively.
The moduli in
\begin{equation*}
    \mathcal{M}_0=\{48, 49,\ldots, 480\}
\end{equation*}
can be used in all $19$ coverings, and each modulus in
\begin{equation*}
    \mathcal{M}_1=\{3, 4,\ldots, 18, 20,\ldots, 22, 24, 25\}
\end{equation*}
can be used in just one of the coverings.\\
\indent Now,
\begin{equation*}
    \sum_{d\in\mathcal{M}_0}\frac{1}{d}\approx0.92669<1,
\end{equation*}
which means that each of the nineteen coverings needs at least one of the $21$ moduli in $\mathcal{M}_1$ to get the sum of the reciprocals to be $1$ or more.\\
\indent After relabeling, we can assume that $3\in\mathcal{M}'_0,\ 4\in\mathcal{M}'_1,\ldots,\ \text{and}\ 13\in\mathcal{M}'_{10}$. This leaves us with ten moduli, $14, 15, 16, 17, 18, 20, 21, 22, 24, 25$, to distribute into the remaining eight coverings, but none of these remaining moduli alone are enough to get the sum of the reciprocals to be $1$ or more, and thus we have a contradiction. Therefore, the theorem holds for $m=48$.\\
\hfill \break
\indent For $m=49$, suppose we have a distinct covering with the moduli of all the congruences in $[49, 490]$. We have,
\begin{equation*}
    \sum_{\substack{d|L_{49} \\ 49\leq d\leq490}}\frac{1}{d}\approx1.0247320.
\end{equation*}
\indent Next, denote the resulting covering by $\mathcal{C}'$ and denote the list of the moduli of the congruences in $\mathcal{C}'$ by $\mathcal{M}'$.\\
\indent Next, we reduce the covering $\mathcal{C}'$ modulo $19$ and obtain $19$ coverings, say $\mathcal{C}'_0, \mathcal{C}'_1,\ldots, \mathcal{C}'_{18}$, whose moduli are $\mathcal{M}'_0, \mathcal{M}'_1,\ldots, \mathcal{M}'_{18}$, respectively.
The moduli in
\begin{equation*}
    \mathcal{M}_0=\{49, 50,\ldots, 490\}
\end{equation*}
can be used in all $19$ coverings, and each modulus in
\begin{equation*}
    \mathcal{M}_1=\{3, 4,\ldots,18, 20,\ldots ,25\}
\end{equation*}
can be used in just one of the coverings.\\
\indent Now,
\begin{equation*}
    \sum_{d\in\mathcal{M}_0}\frac{1}{d}\approx0.907898<1,
\end{equation*}
which means that each of the nineteen coverings needs at least one of the $22$ moduli in $\mathcal{M}_1$ to get the sum of the reciprocals to be $1$ or more.\\
\indent After relabeling, we can assume that $3\in\mathcal{M}'_0,\ 4\in\mathcal{M}'_1,\ldots,\ \text{and}\ 10\in\mathcal{M}'_7$. This leaves us with fourteen moduli, $11,\ldots, 18, 20,\ldots, 25$, to distribute into the remaining eleven coverings, but none of these remaining moduli alone are enough to get the sum of the reciprocals to be $1$ or more, and thus we have a contradiction. Therefore, the theorem holds for $m=49$.\\
\hfill \break
For $m=50$, suppose we have a distinct covering with the moduli of all the congruences in $[50, 500]$. We have,
\begin{equation*}
    \sum_{\substack{d|L_{50} \\ 50\leq d\leq500}}\frac{1}{d}\approx1.0083683.
\end{equation*}
Next, denote the resulting covering by $\mathcal{C}'$ and denote the list of the moduli of the congruences in $\mathcal{C}'$ by $\mathcal{M}'$.\\
Next, we reduce the covering $\mathcal{C}'$ modulo $19$ and obtain $19$ coverings, say $\mathcal{C}'_0, \mathcal{C}'_1,\ldots, \mathcal{C}'_{18}$, whose moduli are $\mathcal{M}'_0, \mathcal{M}'_1,\ldots, \mathcal{M}'_{18}$, respectively.\\
The moduli in
\begin{equation*}
    \mathcal{M}_0=\{50, 51,\ldots, 495\}
\end{equation*}
can be used in all $19$ coverings, and each modulus in
\begin{equation*}
    \mathcal{M}_1=\{3, 4,\ldots, 26\}
\end{equation*}
can be used in just one of the coverings.\\
Now,
\begin{equation*}
    \sum_{d\in\mathcal{M}_0}\frac{1}{d}\approx0.88951<1,
\end{equation*}
which means that each of the nineteen coverings needs at least one of the $23$ moduli in $\mathcal{M}_1$ to get the sum of the reciprocals to be $1$ or more.\\
After relabeling, we can assume that $3\in\mathcal{M}'_0,\ 4\in\mathcal{M}'_1,\ldots,\ \text{and}\ 9\in\mathcal{M}'_6$. This leaves us with fifteen moduli, $10,\ldots, 26$, to distribute into the remaining twelve coverings, but none of these remaining moduli alone are enough to get the sum of the reciprocals to be $1$ or more, and thus we have a contradiction. Therefore, the theorem holds for $m=50$.
\subsection*{Reduction via Lemma \ref{lem3} and Lemma \ref{lem4}}
For all $m\in\{7, 21\}$, we can get the sums of the reciprocals of the moduli to be $<1$ by using both Lemma \ref{lem3} and \ref{lem4}. We write these cases out here.\\
\hfill \break
For $m=7$, suppose we have a distinct covering with the moduli of all the congruences in $[7, 70]$. We have,
\begin{equation*}
    \sum_{\substack{d|L_7 \\ 7\leq d\leq70}}\frac{1}{d}\approx1.2396825.
\end{equation*}
\indent Since there are $2$ permissible multiples of $2^4$ in the interval $[7, 70]$, namely $2^4$ and $2^4\cdot3$, we can replace the congruences with these moduli by a single congruence modulo $2^3\cdot3$ and still have a covering. This brings the sum of the reciprocals down to approximately $1.198016$.\\
\indent Denote the resulting covering by $\mathcal{C}'$ and denote the list of the moduli of the congruences in $\mathcal{C}'$ by $\mathcal{M}'$.\\
\indent Next, we reduce the covering $\mathcal{C}'$ modulo $7$ and obtain $7$ coverings, say $\mathcal{C}'_0, \mathcal{C}'_1,\ldots, \mathcal{C}'_6$, whose moduli are $\mathcal{M}'_0, \mathcal{M}'_1,\ldots, \mathcal{M}'_6$, respectively. The moduli in
\begin{equation*}
    \mathcal{M}_0=[8, 9, 10, 12, 15, 18, 20, 24, 30, 36, 40, 45, 60, 24]
\end{equation*}
can be used in all $7$ coverings, and each modulus in
\begin{equation*}
    \mathcal{M}_1=\{1, 2, 3, 4, 5, 6, 8, 9, 10\}
\end{equation*}
can be used in just one of the coverings.\\
\indent Now,
\begin{equation*}
    \sum_{d\in\mathcal{M}_0}\frac{1}{d}=0.8<1,
\end{equation*}
which means that each of the seven coverings needs at least one of the $9$ moduli in $\mathcal{M}_1$ to get the sum of the reciprocals to be $1$ or more.\\
\indent After relabeling, we can assume that $1\in\mathcal{M}'_0,\ 2\in\mathcal{M}'_1,\ldots,\ \text{and}\ 4\in\mathcal{M}'_3$. This leaves us with five moduli, $5, 6, 8, 9, 10$, to distribute into the remaining three coverings, since these remaining moduli get the sum of the reciprocals to be $1$ or more, we need to reduce again modulo $p$. Since the density of $\mathcal{M}'_4$ is exactly $1$, we reduce $\mathcal{M}'_4$ modulo $5$.\\
\indent Denote the resulting covering by $\mathcal{C}''$ and denote the list of the moduli of the congruences in $\mathcal{C}''$ by $\mathcal{M}''$.\\
Next, we reduce the covering $\mathcal{C}''$ modulo $5$ and obtain $5$ coverings, say $\mathcal{C}''_0, \mathcal{C}''_1,\ldots, \mathcal{C}''_4$, whose moduli are $\mathcal{M}''_0, \mathcal{M}''_1,\ldots, \mathcal{M}''_4$, respectively. The moduli in
\begin{equation*}
    \mathcal{M}_0=[8, 9, 12, 18, 24, 36, 24]
\end{equation*}
can be used in all $5$ coverings, and each modulus in
\begin{equation*}
    \mathcal{M}_1=\{1, 2, 3, 4, 6, 8, 9, 12\}
\end{equation*}
can be used in just one of the coverings.\\
\indent Now,
\begin{equation*}
    \sum_{d\in\mathcal{M}_0}\frac{1}{d}\approx0.486111<1,
\end{equation*}
which means that each of the five coverings needs at least one of the $8$ moduli in $\mathcal{M}_1$ to get the sum of the reciprocals to be $1$ or more.\\
\indent After relabeling, we can assume that $1\in\mathcal{M}''_0$. This leaves us with seven moduli, $2, 3, 4, 6, 8, 9, 12$, to distribute into the remaining four coverings. We distribute $2, 12$ into $\mathcal{M}''_1$ and distribute $3, 9, 8$ into $\mathcal{M}''_2$. This leaves us with two moduli, but none of these remaining moduli alone are enough to get the sum of the reciprocals to be $1$ or more, and thus we have a contradiction. Therefore, the theorem holds for $m=7$.\\
\hfill \break
For $m=21$, suppose we have a distinct covering with the moduli of all the congruences in $[21, 210]$. We have,
\begin{equation*}
    \sum_{\substack{d|L_{21} \\ 21\leq d\leq210}}\frac{1}{d}\approx1.1161306.
\end{equation*}
\indent Since there are $2$ permissible multiples of $2^6$ in the interval $[21, 210]$, namely $2^6$ and $2^6\cdot3$, we can replace the congruences with these moduli by a single congruence modulo $2^5\cdot3$ and still have a covering. This brings the sum of the reciprocals down to approximately $1.156349$.\\
\indent Denote the resulting covering by $\mathcal{C}'$ and denote the list of the moduli of the congruences in $\mathcal{C}'$ by $\mathcal{M}'$.\\
\indent Next, we reduce the covering $\mathcal{C}'$ modulo $13$ and obtain $13$ coverings, say $\mathcal{C}'_0, \mathcal{C}'_1,\ldots, \mathcal{C}'_{12}$, whose moduli are $\mathcal{M}'_0, \mathcal{M}'_1,\ldots, \mathcal{M}'_{12}$, respectively. The moduli in
\[
    \mathcal{M}_0=[21, 22,\ldots, 96,\ldots, 210, 96]
\]
can be used in all $13$ coverings, and each modulus in
\begin{equation*}
    \mathcal{M}_1=\{2, 3, 4, 5, 6, 7, 8, 9, 10, 11, 12, 14, 15, 16\}
\end{equation*}
can be used in just one of the coverings.\\
\indent Now,
\begin{equation*}
    \sum_{d\in\mathcal{M}_0}\frac{1}{d}\approx0.928498<1,
\end{equation*}
which means that each of the thirteen coverings needs at least one of the $14$ moduli in $\mathcal{M}_1$ to get the sum of the reciprocals to be $1$ or more.\\
\indent After relabeling, we can assume that $2\in\mathcal{M}'_0,\ 3\in\mathcal{M}'_1,\ldots,\ \text{and}\ 12\in\mathcal{M}'_{10}$. This leaves us with three moduli, $14, 15, 16$, to distribute into the remaining two coverings, but none of these remaining moduli alone are enough to get the sum of the reciprocals to be $1$ or more, and thus we have a contradiction. Therefore, the theorem holds for $m=21$.\\
\hfill \break
Because the sum of the reciprocals of the moduli are $<1$ for all $m\in\{6, 7, 8, 9, 15, 16, 18, 20, 21, 22, 23, 24, 25, 33, 42, 43, 44, 45, 46, 47, 48, 49, 50\}$, we conclude for each integer $m\geq3$, there is no distinct covering system with all moduli in the interval $[m, 10m]$.
\end{proof}

\section{Open problems and further work}
As already mentioned the study of the nonexistence of distinct coverings with all moduli in the interval $[n, kn]$ for an integer $k>1$ and $n\geq3$ has a rich history. Beginning in 2006, M. Filaseta, O. Trifonov, and G. Yu answered the question with all moduli in the interval $[n, 6n]$ for $n\geq3$ \cite{Filaseta2}. This was improved to $[n, 8n]$ by the first author and O. Trifonov, \cite{Dalton2}, followed by the first author proving it for $[n, 9n]$ \cite{Dalton1}. Theorem \ref{thm2} proves it for $[n, 10n]$ for $n\geq3$. The natural question therefore is whether there exists a finite upper bound to which there \textit{does} exist a distinct covering with all moduli in the interval $[n, kn]$. This question was answered in the negative by M. Filaseta, K. Ford, S. Konyagin, C. Pomerance, and G. Yu in 2007 \cite{Filaseta1}. In particular, they showed the following theorem,
\begin{thm}[Filaseta et. al., 2007]\label{thm3}
    For any number $K>1$ and $N$ sufficiently large, depending on $K$, there is no covering system using distinct moduli from the interval $(N, KN]$.
\end{thm}
Theorem \ref{thm3} is now vacuously true because of the breakthrough by R. Hough in 2015 \cite{Hough}. However, the study of intervals for the nonexistence of covering systems with distinct moduli is still important. In particular, in private communication with M. Filaseta it would be nice to show there does not exist a distinct covering with all moduli in the interval $[n, 12n]$ for $n\geq4$.\\
\indent The next interesting interval would be to show there exists no distinct covering with all moduli in the interval $[n, 15n]$ for $n\geq5$ but we leave this for further work.\\
\indent Lastly, if one would want to extend our work to larger intervals such as $[n, 12n]$ or $[n, 15n]$ the number of values for which the sum of the reciprocals of the moduli is $>1$ would be quite larger. A suitable modification of this method may suffice, as this current method is undoubtedly not optimal. Perhaps, one may be able to improve on our work using the well-known distortion method of Balister et al. If you are curious, we refer you to an expository paper on the distortion method here \cite{Balister2}. We believe that it would be worthwhile to explore other avenues such as the distortion method or a suitable modification of this method.

\subsection*{Acknowledgment}  The second author would like to thank M. Filaseta and O. Trifonov for helpful conversations pertaining to this research during a conference at the University of South Carolina.

\bibliographystyle{plain}

\begin{thebibliography}{99}
		
	
	\bibitem{Balister1}
	P. Balister, B. Bollob\'as, R. Morris, J. Sahasrabudhe and M. Tiba,\,
	{\it On the Erd\H{o}s covering problem : The density of the uncovered set.} 
	Invent. Math. 228 (2022), no. 1, pp. 377--414.


    \bibitem{Balister2}
	P. Balister, B. Bollob\'as, R. Morris, J. Sahasrabudhe and M. Tiba,\,
	{\it Erd\H{o}s covering systems.} 
	Acta. Math. Hungar. 161 (2020), pp. 540--549.
    
    
    \bibitem{Erdos}
    P. Erd\H os,\,
    {\it On integers of the form $2^k+p$ and some related problems.} 
    Summa Brasil.~Math.~2 (1950), pp. 113--123.


    \bibitem{Dalton1}
    J. Dalton,\,
    {\it Extreme Covering Systems. Primes Plus Squarefrees, and Lattice Points Close to a Helix.} 
    Thesis (Ph.D.)--University of South Carolina. ProQuest LLC, Ann Arbor, MI, 2023.


    \bibitem{Dalton2}
    J. Dalton and O. Trifonov,\,
    {\it Extreme Covering Systems.} 
    J. Integer Seq.~(9) 25 (2022), pp. 1--28.


    \bibitem{Filaseta1}
    M. Filaseta, K. Ford, S. Konyagin, C. Pomerance and G. Yu,\,
    {\it Sieving by large integers and covering systems of congruences.} 
    J. Amer. Math. Soc.~20 (2007), pp. 495--517.


    \bibitem{Filaseta2}
    M. Filaseta, O. Trifonov, G. Yu,\,
    {\it Distinct covering systems,} 
    research seminar, 2006.


    \bibitem{Krukenberg}
    C. Krukenberg,\,
    {\it COVERING SETS OF THE INTEGERS.} 
    Thesis (Ph.D.)--University of Illinois at Urbana-Champaign. ProQuest LLC, Ann Arbor, MI, 1971.
    
    
    \bibitem{Hough} 
    R.~Hough,\,
    {\it Solution of the minimum modulus problem for covering systems.}
    Ann.~of Math.~(2) 181 (2015),\, no.~1, pp. 361--382.

\end{thebibliography}

\end{document}